\documentclass[11pt]{article}

\usepackage[utf8]{inputenc}
\usepackage[T1]{fontenc}
\usepackage{times}
\usepackage{amsmath,amsfonts,amssymb}
\usepackage{algorithm}
\usepackage{algorithmic}
\usepackage{empheq}
\usepackage{graphicx}
\usepackage{caption}
\usepackage{subcaption}
\usepackage{booktabs}
\usepackage{xcolor}
\usepackage{hyperref}
\usepackage{enumitem}
\usepackage{cite}

\usepackage[margin=1in]{geometry}

\begin{document}

\title{An Additively Preconditioned Trust Region Strategy for Machine Learning}

\author{
Samuel Cruz Alegría\thanks{Università della Svizzera Italiana} \quad
Bindi Çapriqi\thanks{King Abdullah University of Science and Technology. Corresponding author: \texttt{bindi.capriqi@kaust.edu.sa}} \quad
Shega Likaj\textsuperscript{1,2} \\
Ken Trotti\textsuperscript{2} \quad
Rolf Krause\textsuperscript{1,2} \\[1em]
\small \textsuperscript{1}Università della Svizzera Italiana \\
\small \textsuperscript{2}King Abdullah University of Science and Technology
}

\date{}

\maketitle

\begin{abstract}
Modern machine learning, especially the training of deep neural networks, depends on solving large-scale, highly nonconvex optimization problems, whose objective function exhibit a rough landscape. 
Motivated by the success of parallel preconditioners in the context of Krylov methods for large scale linear systems, we introduce a 
novel nonlinearly preconditioned Trust-Region method that makes use of an additive Schwarz correction at each minimization step,
thereby accelerating convergence.

More precisely, we propose a variant of the Additively Preconditioned Trust-Region Strategy (APTS), which combines a right-preconditioned additive Schwarz framework with a classical Trust-Region algorithm. By decomposing the parameter space into sub-domains, APTS solves local non-linear sub-problems in parallel and assembles their corrections additively. The resulting  method not only shows fast convergence; due to the underlying Trust-Region strategy, it 
furthermore largely obviates the need for hyperparameter tuning.
\end{abstract}

\vspace{1em}
\noindent

\section{Introduction}

Given a dataset \(\{(x^{(i)},y^{(i)})\}_{i=1}^M\), with $x^{(i)} \in \mathbb{R}^q$ and $y^{(i)} \in \mathbb{R}^p$, in supervised learning  parameters $\theta \in \mathbb{R}^n$ are to be identified such that the following minimization problem is approximately solved:
\begin{equation}\label{def:min-problem}
\min_{\theta \in \mathbb{R}^n} f(\theta) 
\ :=\ 
\min_{\theta \in \mathbb{R}^n} 
 \frac{1}{M} \sum_{i=1}^M \ell\bigl(\mathcal{N}({\theta};x^{(i)}), y^{(i)}\bigr). 
\end{equation} 
Here $f: \mathbb{R}^n \longrightarrow \mathbb{R}$ is the empirical risk or average loss over the whole dataset, \(\ell:\mathbb{R}^p\times\mathbb{R}^p\longrightarrow\mathbb{R}\) denotes the per‐sample loss and \(\mathcal{N}({\theta}, \cdot):\mathbb{R}^q\longrightarrow\mathbb{R}^p\) is the forward pass of a Neural Network (NN), that maps an input vector $x\in \mathbb{R}^q$ to a prediction $y\in \mathbb{R}^p$ under the model’s parameters $\theta$. 
In the context of NNs, approximately solving (\ref{def:min-problem}) is known as training.

Loss functions in Machine Learning are typically large scale, i.e. defined on large data-sets, and are highly nonconvex  with varying curvature across different layers or parameter blocks.
They often  admit for many saddle points and poor local minima, i.e. minima with insufficient approximation properties of the resulting network. Moreover, training usually starts from random or heuristic initializations, which can lie far from a ``good'' solution.  Under these conditions, naïve gradient updates may stall or converge very slowly, on extended plateaus around high-dimensional saddle points \cite{Dauphin-et-al-2014,Jin-et-al-2017}. Therefore, ensuring stable, hyperparameter‐robust descent from arbitrary starting points remains a central challenge for scalable ML optimization and motivates our exploration of additive preconditioned Trust Region (TR) methods.

The most widely used optimization methods for problem \eqref{def:min-problem} are Stochastic Gradient Descent (SGD)  and its variants\cite{Zinkevich-2010,Bottou-Curtis-Nocedal-2018,Dai-Zhu-2018}, most notably Adaptive Moment estimation (Adam) \cite{Kingma-Ba-2014}. Further, the limited-memory Broyden–Fletcher–Goldfarb–Shanno (L‑BFGS) quasi‑Newton method is particularly effective for ill‑conditioned, large‑scale problems, since it builds an efficient low‑rank approximation of the inverse Hessian matrix \cite{Nocedal-Wright-1999}.

Although  these methods are widely used and have wide-ranging benefits,  they also have their specific drawbacks. For SGD, each iteration is inexpensive, but an asymptotically diminishing 
learning‑rate schedule is required to guarantee convergence even in convex settings \cite{Bottou-Curtis-Nocedal-2018}. Variance-reduced variants such as \ SVRG \cite{Johnson-Zhang-2013}, SAGA \cite{Defazio-Bach-Lacoste-Julien-2014}, and SARAH \cite{Nguyen-Liu-Scheinberg-Takac-2017} alleviate this by controlling the 
gradient noise and can attain linear convergence without 
vanishing stepsizes for finite-sum objectives, yet they still 
incur periodic full-gradient costs and their gains for deep, highly nonconvex networks are mixed \cite{Defazio-Bottou-2019}.  
Adaptive schemes such as Adam remove the need for manual decay, 
yet still lack global convergence guarantees for nonconvex 
settings, and may face issues even in convex settings 
\cite{Kingma-Ba-2014,Reddi-Kale-Kumar-2019,Ruder-2016}. Momentum-accelerated SGD variants are ubiquitous but likewise sensitive to 
stepsize choice \cite{Sutskever-Martens-Dahl-Hinton-2013}. On the 
other hand, quasi-Newton methods like L‑BFGS are better‑suited 
than pure first‑order methods for ill‑conditioned problems, but 
more expensive than SGD variants and require non‑trivial 
adaptations for the stochastic setting  \cite{Berahas-Nocedal-Takac-2016, Berahas-Jahani-Richtarik-Takac-2022}.

Trust-Region (TR) methods are globalization strategies that improve robustness by adapting the TR radius dynamically and offer global convergence properties for both convex and nonconvex settings \cite{Nocedal-Wright-1999, connGouldToint2000}. They are less common in the ML context, but are becoming increasingly popular. However, this class of methods also requires non-trivial adaptations to be applied to the multi-level \cite{Gross-Krause-2009} and/or stochastic settings \cite{Rafati-DeGuchy-Marcia-2018, Curtis-Scheinberg-Shi-2019, Gao-Ng-2022, Gratton-Jerad-Toint-2023, Gratton-Kopanicakova-Toint-2023, Kopanicakova-Krause-2022}.

Another important issue in this context is how to efficiently parallelize the training process, which becomes increasingly challenging with growing data and network complexity \cite{Catlett-1991, Blum-Rivest-1988, Srivastava-Greff-Schmidhuber-2015, Grosse-2018}. Many modern learning problems involve extremely high-dimensional parameter spaces.  When training such models, a global update can become expensive both computationally and memory-wise. There are two main approaches to parallelization. One is the data parallel approach, which can be described as splitting the dataset into chunks that can be used to train model copies in parallel. This approach has been applied to various methods \cite{CruzAlegría-et-al-2025} and has proven useful for very large datasets \cite{Shallue-2018}. The other main approach is the model parallel approach, which involves splitting the model into separate parts (e.g. each layer) which can be trained in parallel \cite{BenNun-Hoefler-2019}.

Domain decomposition (DD) is a promising approach for addressing parallelization \cite{Shallue-2018, Nichols-Singh-Lin-Bhatele-2021, BenNun-Hoefler-2019, Trotti-et-al-2025}. In the context of model-parallelism, the techniques can split the overall parameter vector into smaller pieces that can be optimized independently and in parallel. DD methods were originally developed for solving large-scale partial differential equations through decomposition of the global problem into subproblems, which can be solved independently and used as aid for the main problem \cite{Toselli-Widlund-2004}. The additive approach to DD allows the subproblems to be solved independently and thus naturally leads to parallelization \cite{Chan-Zou-1994, Toselli-Widlund-2004, Gander-2006, Erhel-Gander-Halpern-Pichot-Sassi-Widlund-2014}. 

Notable variants include Additive Schwarz Preconditioned Inexact Newton (ASPIN) \cite{Cai-Keyes-2002}, Globalized ASPIN (GASPIN) \cite{Gross-Krause-2021}, Restricted Additive Schwarz Preconditioned Exact Newton (RASPEN) \cite{Dolean-Gander-Kheriji-Kwok-Masson-2016}, and the Additively Preconditioned Trust-Region Strategy (APTS) \cite{Gross-Thesis,CruzAlegría-et-al-2025,Trotti-et-al-2025}. Most of the additive techniques mentioned were designed to solve nonlinear systems of equations; however, APTS and GASPIN were specifically designed for nonconvex problems, and both provide global convergence guarantees through a TR globalization strategy. DD Methods have been applied to a much wider variety of problems in the ML context. For a comprehensive survey, consult \cite{Klawonn-Lanser-Weber-2023}.

In this paper, we present a novel and cheaper adaptation of APTS tailored for ML applications \cite{Gross-Thesis,CruzAlegría-et-al-2025, Trotti-et-al-2025}. Our approach addresses three critical challenges that hinder many popular ML optimizers: providing global convergence guarantees in deterministic settings, reducing reliance on laborious hyperparameter tuning, and enabling efficient parallelization. Through this framework, we aim to combine the robustness of TR methods with the scalability demands of modern NN training.

The remainder of the paper is organized as follows. In Section \ref{sec:theory} we introduce the theoretical framework, covering nonlinear right preconditioning and DD principles.  Section \ref{sec:apts} details the APTS method, including the global preconditioner and TR loop.  Section \ref{subsec:implementation} discusses implementation details such as approximated subdomains and GPU‐parallel execution. Numerical experiments on well-established datasets like the MNIST dataset and a convolutional NN (CNN) are presented in Section \ref{sec:results}.  Finally, Section \ref{sec:conclusion} presents the conclusions and outlines future work.

\section{Theoretical Framework of APTS}
\label{sec:theory}

In this section, we describe the theoretical foundations of APTS. We first reformulate the global minimization problem by means of nonlinear right preconditioning. Next, we transpose classical DD principles to the parameter space of deep NNs. Finally, we combine these elements into an additive preconditioner that operates within a TR globalization framework. This section introduces the notation, assumptions, and auxiliary results required to motivate the algorithmic design presented in Section \ref{sec:apts}.

\subsection{Nonlinear Right Preconditioning}
\label{subsection:nonlinear-right-preconditioning}

When using iterative methods to solve linear systems, given an iterate $\theta \in \mathbb{R}^n$, at each iteration we seek a step \(s\in\mathbb R^{n}\) that approximately solves
\[
A\,(\theta+s) \;=\; b. 
\]
In many large‐scale linear problems, the matrix \(A\) exhibits poor spectral properties, leading to slow convergence and numerical instability.  A classical remedy is to introduce a right preconditioner \(M\), effectively introducing the change of variable \(\theta +s = M^{-1}\,\tilde{\theta}\). Here, the preconditioner \(M\) is chosen so that the spectrum of \(A\,M\) is significantly more clustered than that of \(A\) alone, thereby accelerating the convergence of the inner linear solve and producing a more robust update \(\tilde \theta\).  In this framework, for a given iterate \(\theta \in \mathbb{R}^n\) we seek an incremental update \(s \in \mathbb{R}^n\) by solving

\begin{subequations}\label{eq:right-preconditioning}
\begin{empheq}[left=\empheqlbrace]{align}
A\,M\bigl(\theta + s\bigr) \;-\; b &= 0, \\
\tilde{\theta} &= M\bigl(\theta + s\bigr).
\end{empheq}
\end{subequations}

The same idea can be applied to nonlinear optimization. Unconstrained optimization can be cast as a nonlinear root–finding problem: seeking a critical point \(\theta^\star\) of a smooth objective  \(f\colon\mathbb{R}^n\!\longrightarrow\mathbb{R}\) is equivalent to solving the  nonlinear system
\[
       \nabla f(\theta)\;=\; 0.
\]
Analogous to right preconditioning in linear systems, nonlinear right preconditioning transforms the problem by applying a nonlinear map to the argument of the gradient. Specifically, instead of solving \(\nabla f(\theta) = 0\) directly, we introduce an auxiliary operator
\[
  \mathcal{F}:\mathbb{R}^n\longrightarrow\mathbb{R}^n, \quad \mathcal{F}\in C^1(\mathbb{R}^n),
\]
which is called \emph{nonlinear preconditioner}, and reformulate the problem as
\begin{subequations}\label{eq:nonlinear_precond_operator}
\begin{empheq}[left=\empheqlbrace]{align}
  \nabla f(\mathcal{F}(\theta+s)) &= 0, \label{eq:nonlinear_precond_a}\\
  \theta &= \mathcal{F}(\theta+s). \label{eq:nonlinear_precond_b}
\end{empheq}
\end{subequations}
We now apply Newton's method \cite{Nocedal-Wright-1999} to equation \eqref{eq:nonlinear_precond_a} to find \(s\), obtaining:
\begin{subequations}\label{eq:nonlinear_precond}
\begin{empheq}[left=\empheqlbrace]{align}
    \nabla^2 f(\mathcal{F}\theta^k))\,\mathcal{F}'(\theta^k)\,s &= -\nabla f(\mathcal{F}(\theta^k)), \label{eq:nonlinear_newton_a}\\
    \theta^{k+1} &= \mathcal{F}(\theta^k+s), \label{eq:nonlinear_newton_b}
\end{empheq}
\end{subequations}
where \(\theta^k\) is the current iterate and \(s\) is the search direction and \(\mathcal{F}'\) is the Jacobian of \(\mathcal{F}\).\\
The evaluation of \(\mathcal{F}(x^k+s)\) can be computationally expensive, especially if \(\mathcal{F}\) is a complex nonlinear function. To mitigate this, we can linearize \cite{Gross-Thesis} the update step by approximating equation \eqref{eq:nonlinear_newton_b} as
\[
    \mathcal{F}(\theta^k+s) \approx \mathcal{F}(\theta^k) + s^k.
\]
with \(s^k:=\mathcal{F}'(\theta^k)\!\,s\), thus obtaining
\begin{subequations}\label{eq:linearized_nonlinear_precond}
\begin{empheq}[left=\empheqlbrace]{align}
    \nabla^2 f(\mathcal{F}(\theta^k))\!\,s^k &= -\nabla f(\mathcal{F}(\theta^k)), \label{eq:linearized_nonlinear_newton_a}\\
    \theta^{k+1} &= \mathcal{F}(\theta^k) + s^k, \label{eq:linearized_nonlinear_newton_b}
\end{empheq}
\end{subequations}
where \(s^k\) is obtained through the application of the TR method.

\subsection{Parameter Decomposition}
To adapt the DD strategy from PDEs to NN training, we partition the parameter index set 

\[
   \{1,\dots,n\} = \bigcup_{d=1}^N C_d,
\]
so that each subset \(C_d\) has size \(\lvert C_d\rvert=n_d\). In this paper, we will only consider the non-overlapping case, i.e., \(C_d\cap C_{d'}=\emptyset\) for \(d\neq d'\), noting that overlapping partitions can be considered and will only lead to minor modifications to the transfer operators. Each subset $C_d$ defines the subdomain $ D_d =  \mathbb{R}^{n_d} \subset \mathbb{R}^n$.
Let us define the transfer operators, specifically the restriction and prolongation maps
\[
  R_d:\mathbb{R}^n\longrightarrow D_d, 
  \qquad
  R_d^T:D_d\longrightarrow\mathbb{R}^n.
\]
Here \(R_d\) downsamples the parameter vector \(\theta\in\mathbb{R}^n\) to the components in \(C_d\), enabling the definition of the local subproblem on each block, and \(R_d^T\) upsamples it back to \(\mathbb{R}^n\) by zero‑padding the components not in \(C_d\), which allows the independent updates to be recombined into the full parameter vector. We note that, by construction, 
\[
  R_d R_d^T = I_{n_d},
  \qquad
  R_d R_{d'}^T = 0_n\;\;(d\neq d'),
  \qquad \text{and}\qquad
  \sum_{d=1}^N R_d^T R_d = I_n,
\]
where \(I_n\) and $0_n$ are the identity and null operators on \(\mathbb{R}^n\), respectively.

\subsection{Building the Additive Preconditioner}

Through the transfer operators, we can now define the restricted objective function $f(\theta^k)$ on a subdomain, i.e.,
\begin{equation}\label{eq:local_function}
  f_d(\theta_d) = f(R_d^T\theta_d + (I_n - R_d^T R_d)\theta^k),\qquad\theta^k \in\mathbb{R}^n
\end{equation}
which consists in function $f$ with variables $(\theta^k)_i$ for $i \notin C_d$ that become frozen parameters.

Following the approach originally proposed in \cite{Nash-2000}, in order to maintain the first-order consistency with the original minimization problem in equation \eqref{def:min-problem} we have to add the first-order consistency term which consists of the difference between the restricted gradient of the global objective evaluated at the global iterate $R_d\nabla f(\theta^k)$ and the gradient of the restricted objective evaluated at the restricted global iterate $\nabla f_d(R_d\theta^k)$. Thus, we define the local objective as

\begin{equation} \label{eq:local_objective}
  \tilde f_d(\theta_d)
  = \underbrace{f_d\;\bigl(\theta_d)}_{f \  \text{restricted to subdomain} \ D_d}
  + \underbrace{\langle R_d\nabla f(\theta^k)-
\nabla f_d(R_d\theta^k),\theta_d-R_d\theta^k\rangle}_{\text{first-order correction term}},
\end{equation}

Now, by applying the preconditioning strategy developed in Section \ref{subsection:nonlinear-right-preconditioning}, we define the subdomain preconditioner as $\mathcal{F}_d:D_d\rightarrow D_d$ and the subdomain step obtained after the application of $\mathcal{F}_d$ as $s_d^k:=\mathcal{F}_d(R_d\theta^k)-R_d\theta^k$.

The global preconditioner $\mathcal{F}$ consists of the parallel and decoupled execution of the local preconditioners $\mathcal{F}_d,\ i=1,...,N$ and a strategy to combine the steps. Formally:

\begin{equation} \label{def:precondition-operator}
\mathcal{F}\bigl(\theta^k\bigr)
\;=\;
\mathcal A\Bigl( R_1^\mathrm{T}s_1^k,\;\dots,\; R_Ns_N^k,\; \theta^k \Bigr),
\end{equation}

where $\mathcal A$ is the recombination strategy which, in the case of the APTS algorithm, is the following
\[
\mathcal A_{\mathrm{APTS}}\bigl(R_1^T s_1^k,\dots,R_N^T s_N^k,\theta^k\bigr)
=
\begin{cases}
\displaystyle
\theta^k \;+\;\sum_{d=1}^N R_d^T s_d^k,
&\text{if }\sum_{d=1}^N R_d^T s_d^k \;\;\;\text{is “sufficiently good”},\\[1ex]
\theta^k,
&\text{otherwise}.
\end{cases}
\]

\color{black}

\subsection{The TR Method}\label{subsec:TR}

\color{black}
The TR method, summarized in Algorithm~\ref{alg:tr}, is an iterative algorithm that, at each iteration, computes a candidate step and then decides whether to accept it. To do so, the method replaces the full objective $f$ by a local model $m^k$, which is assumed to be a reliable approximation of $f$ only in a neighborhood of the current iterate $k$. This neighborhood, known as the trust region, is defined as a Euclidean ball:
\[
  \mathcal B^k \;:=\; \bigl\{\,s\in\mathbb{R}^n : \|s\|_2 \le \Delta^k \bigr\},
  \quad \Delta^k>0,
\]
where \(\Delta^k\) is the TR radius and \(s\) the candidate update step. We adopt the standard quadratic model, i.e., the second-order Taylor approximation of $f$:

\begin{equation}\label{eq:TR-model}
  m^k(s)\;=\;
  \nabla f(\theta^k)^\top s
  \;+\;
  \tfrac12\,s^\top H^k s,
  \quad H^k\simeq\nabla^2 f(\theta^k),
\end{equation}
with \(H^k\in\mathbb{R}^{n\times n}\) symmetric.  In practice, \(H^k\) may be the exact Hessian, a quasi‑Newton approximation, or a subsampled Hessian (see e.g. BFGS, SR1)~\cite{Nocedal-Wright-1999}. The step $s^k$ is obtained by minimizing $m^k(s)$ over \(\mathcal B^k\). Once $s^k$ is computed, the quality of the model prediction is assessed via the agreement ratio
\begin{equation}\label{eq:rho}
  \rho^k
  \;=\;
  \frac{f(\theta^k)-f(\theta^k+s^k)}
       {m^k(0)-m^k(s^k)}\in\mathbb{R},
\end{equation}
where $f(\theta^k)-f(\theta^k+s^k)$ is the \emph{actual decrease} in the
objective and $m_k(0)-m_k(s^k)$ is the \emph{predicted decrease}
given by the quadratic model. Thus $\rho^k$ quantifies the accuracy of the quadratic model: \(\rho^k\approx1\) indicates a near–perfect prediction, \(\rho^k\ll1\) signals an over-optimistic model, and \(\rho^k<0\) means the trial step actually increased the objective value. We then compare \(\rho^k\) to two thresholds, \(0<\eta_1<\eta_2<1\), and adjust the iterate and TR radius \(\Delta^k\) using increase/decrease factors \(0<\gamma_{\mathrm{dec}}<1<\gamma_{\mathrm{inc}}\): if \(\rho^k\ge\eta_2\) we accept \(s^k\) and set \(\Delta^{k+1}=\gamma_{\mathrm{inc}}\,\Delta^k\); if \(\eta_1<\rho^k<\eta_2\) we accept \(s^k\) but leave \(\Delta^k\) unchanged; otherwise we reject \(s^k\) and shrink the region via \(\Delta^{k+1}=\gamma_{\mathrm{dec}}\,\Delta^k\).

\color{black}
\begin{algorithm}[htbp]
\caption{\textbf{TR} Trust–Region Method}
\label{alg:tr}
\begin{algorithmic}[1]
  \REQUIRE objective $f(\theta)$, initial iterate $\theta^0$, 
           initial radius $\Delta^0>0$, acceptance threshold $0<\eta_1<\eta_2<1$, 
           decrease and decrease factors $0<\gamma_{\text{dec}}<1<\gamma_{\text{inc}}$, norm $p$, and iteration count $m$.
  \FOR{$k=0,1,2,\dots,m$}
      \STATE $m_k(s) = \nabla f(\theta^k)^\top s + \tfrac{1}{2} s^\top H_k s$, \quad $H_k \simeq \nabla^2 f(\theta^k)$\quad as in equation \eqref{eq:TR-model}.
      \STATE $s^k = \underset{{\|s\|_p\le\Delta^k}}{\arg \min} \; m_k(s)$.

    \IF{$\rho^k\ge\eta_2$} \label{line:accept-start}
      \STATE $\theta^{k+1}\gets\theta^k+s^k,\quad\Delta^{k+1}\gets\gamma_{\mathrm{inc}}\Delta^k$
    \ELSIF{$\eta_1 <  \rho^k<\eta_2$}
      \STATE $\theta^{k+1}\gets\theta^k+s^k,\quad\Delta^{k+1}\gets\Delta^k$
    \ELSE
      \STATE $\theta^{k+1}\gets\theta^k,\quad\Delta^{k+1}\gets\gamma_{\mathrm{dec}}\Delta^k$ 
    \ENDIF  \label{line:accept-end}
  \ENDFOR
\end{algorithmic}
\end{algorithm}

\color{black}

\section{The APTS Method}
\label{sec:apts}

The APTS method, shown in Algorithm ~\ref{alg:apts}, combines parallel TR updates in each subdomain with a single global TR safeguard. During each APTS iteration $k$, we progress through four phases: preconditioning, additive recombination, global acceptance, and global TR iterations.

In the first phase (lines~\ref{line:local-for}-\ref{line:local-end}), we perform \(m\) TR iterations on each subdomain \(D_d\) using the first-order consistent model \(\tilde f_d\) defined in equation \eqref{eq:local_objective}. To guarantee that no sequence of local steps exceeds the global TR radius \(\Delta^k_G\), we initialize each subdomain's TR radius at \(\frac{\Delta^k_G}{m}\) and set the increase factor \(\gamma_{\mathrm{inc}}=1\). Even in the extreme case where all \(m\) inner steps take the full length \(\frac{\Delta^k_G}{m}\) in the same direction, their sum $m\cdot\frac{\Delta^k_G}{m} \;=\; \Delta_G^k$ remains on the boundary of the global trust region. Moreover, because each subdomain solve depends only on its restricted subproblem, the local steps can be computed in parallel without requiring cross-domain communication.

The second phase is described in lines \ref{line:local-end} and \ref{line:additive}. We compute the step on each subdomain $D_d$ as the difference between the initial local variable $\theta_d^{k,0}$ and the final local iterate $\theta_d^{k,m}$
   \begin{equation} \label{de:local_s}
           s_d^k \;=\;\theta_d^{k,m} - \theta_d^{k,0},
    \qquad
    d=1,\dots,N.
  \end{equation}
  Then we gather all the steps and form the tentative global step
  \begin{equation} \label{def:global_s}
    s^k \;=\;\sum_{d=1}^N R_d^T\,s_d^k.
  \end{equation}
Note that the global update $\mathcal{F}(\theta^k) = \theta^k +s^k$ corresponds to the action of the right-preconditioning operator defined in equation \eqref{def:precondition-operator}. 

The third phase is described in lines \ref{line:rho} and \ref{line:tr-agreement-ratio}. We compute the agreement ratio $\rho^k$ in equation \eqref{eq:rho} as

  \begin{equation}\label{eq:apts-rho}
    \rho^k
    \;=\;
    \frac{f(\theta^k) - f(\theta^k + s^k)}
         {\displaystyle\sum_{d=1}^N\Bigl[\tilde f_d\bigl(R_d\theta^k\bigr)
                                       - \tilde f_d\bigl(\theta_d^{k,m}\bigr)\Bigr]},
  \end{equation}
  and adapt the TR radius and accept/reject the step accordingly as explained at the end of Section \ref{subsec:TR}.

  The final phase is described in line \ref{line:global-tr}. Here we perform \(m_G\) TR steps on the global loss function \(f\) to handle any residual coupling.

\begin{algorithm}[htbp]
\caption{\textbf{APTS} Nonlinear Additively Preconditioned Trust–Region Strategy}
\label{alg:apts}
\begin{algorithmic}[1]
  \REQUIRE Loss function $f(\theta)$,
            initial guess $\theta^0$,
            inner–loop iters $m$,
            initial TR radius $\Delta^0_G$,
            usual TR parameters $\eta_1,\eta_2$,
            $\gamma_{\mathrm{dec}},\gamma_{\mathrm{inc}}$,
            the number of global TR iterations $m_G$.
  \FOR{$k=0,1,2,\dots$}
    \FOR{subdomains $d=1,\dots,N$ {\bf in parallel}}
    \label{line:local-for}
      \STATE $\theta_{d}^{k,0}\gets R_d\,\theta^k$
      \label{line:local-start}
        \STATE $\theta_{d}^{k,m}=\;$TR($\tilde{f}_d, \theta^{k,0}_d , \frac{\Delta^k_G}{m}, \eta_1, \eta_2, \gamma_{\text{dec}}, 1, m$)  \label{line:local-end}
      \STATE $s_d^k\gets \theta_{d}^{k,m}-\theta_{d}^{k,0}$ \label{line:additive}
    \ENDFOR

    \STATE Set $s^k$ as defined in \eqref{def:global_s} and compute $\rho^k$ as defined in \eqref{eq:apts-rho}.  \label{line:rho}
    \STATE $\theta^{k+\frac{1}{2}}$ and $\Delta_G^{k+\frac{1}{2}}$: Update $\theta^k$ and $\Delta^k$ according to the strategy in lines \ref{line:accept-start} to \ref{line:accept-end} of the TR algorithm \ref{alg:tr} \label{line:tr-agreement-ratio}.    
    \STATE TR($f, \theta^{k+\frac{1}{2}}, \Delta_G^{k+\frac{1}{2}}, \eta_1, \eta_2, \gamma_{\text{dec}}, \gamma_{\text{inc}}, m_G$) \label{line:global-tr}
  \ENDFOR
\end{algorithmic}
\end{algorithm}

In equation \eqref{eq:local_objective} we note that $\tilde f_d(\theta_d)$ is defined through $f_d(\theta_d)$ in equation \eqref{eq:local_function} and $f_d(\theta_d)$ consists in a NN copy with parameters whose index does not belong to $C_d$ which are frozen. In this case, the computation of the gradient still requires a full forward and a potentially full backward pass, thereby incurring prohibitive costs on deep architectures. In particular, when subdomains lie close to the NN input, propagating long Jacobian chains at every subdomain iteration leads to excessive computation time and memory usage \cite{Trotti-et-al-2025}. We therefore develop an inexact version of the APTS algorithm, explained in Section~\ref{subsec:implementation}. 

\subsection{Inexact APTS}
\label{subsec:implementation}
\begin{figure}[htbp]
    \centering
    \includegraphics[width=0.5\linewidth]{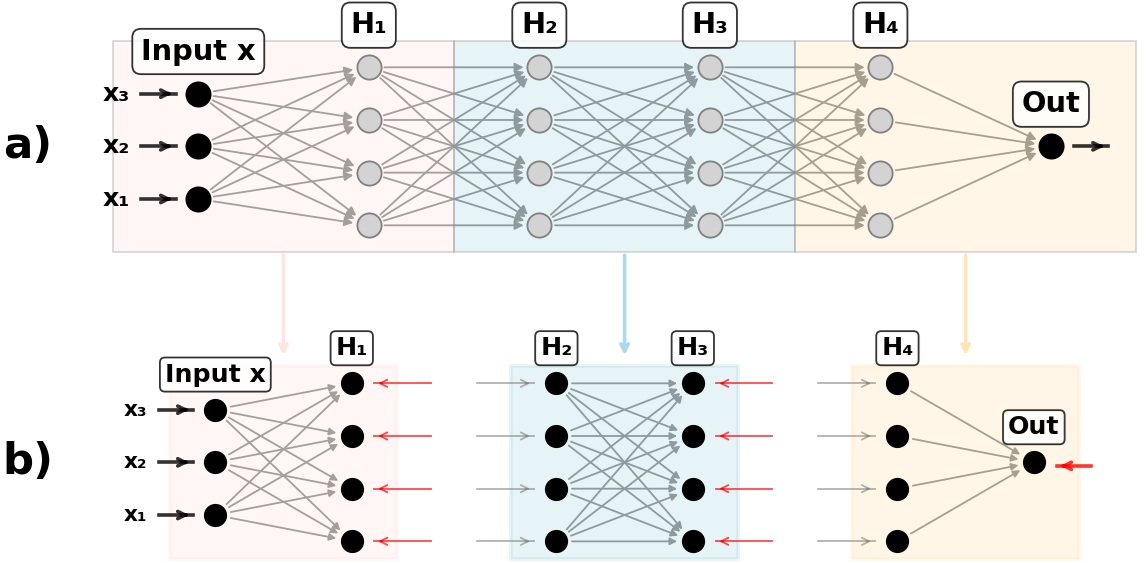}
    \caption{Graphical representation of the pipelined NN across 3 GPUS (a), and the decoupled subdomains (b).}
    \label{fig:nn}
\end{figure}

Rather than defining a subdomain as a full model replica with frozen layers, as $\tilde{f}_d$ in equation~\eqref{eq:local_objective}, 
in the Inexact APTS (IAPTS) variant, we define it structurally as a partition of the network assigned to a specific GPU or group of GPUs. Figure~\ref{fig:nn} illustrates this decomposition in the case of a fully connected network, where colored segments denote the GPU-specific partitions.
IAPTS starts with one global forward and backward before the preconditioning step, during which it stores, for each subdomain, the forward activations coming from the previous subdomain (black arrows in Figure \ref{fig:nn}.b) and the backward gradients received from the subsequent subdomain (red arrows in Figure \ref{fig:nn}.b). 

Therefore, IAPTS caches activations and backward gradients from a single global pass and uses these quantities to reconstruct the required gradient factors locally without reevaluating the frozen layers. For example, let the layers be
\begin{align}
\text{Hidden layer 1:} \quad &H_1(x) = \sigma_1(W_{H_1} x + b_{H_1}), \quad x \text{ is the network input},\\
\text{Hidden layer 2:} \quad &H_2(z) = \sigma_2(W_{H_2} z + b_{H_2}), \\
\text{Hidden layer 3:} \quad &H_3(z) = \sigma_3(W_{H_3} z + b_{H_3}), \\
\text{Hidden layer 4:} \quad &H_4(z) = \sigma_4(W_{H_4} z + b_{H_4}), \\
\text{Output layer:} \quad &\text{Out}(z) = \sigma_5(W_{\text{output}} z + b_{\text{output}}),
\end{align}
where $\sigma_i$ are the activation functions, $W_i$ are weight matrices, and $b_i$ are bias vectors. The full NN is described by
\[
\mathcal{N}(\theta;x) = \text{Out}(H_4(H_3(H_2(H_1(x))))),
\]
with $\theta$ being the vector containing all of the parameters $W_i,b_i$ of the NN.

For a loss $f(y,\hat y)$ with $\hat y=\mathcal{N}(\theta;x)$, gradients follow from the chain rule. We factor each gradient into a stored downstream gradient, obtained during the global pass and frozen locally, multiplied by a local derivative:
\[
\begin{array}{c|c|c}
\text{Parameter} & \text{Stored downstream gradient} & \text{Gradient} \\ \hline
W_{\text{output}} &
G_{\text{out}} = \dfrac{\partial f}{\partial \hat y} &
\dfrac{\partial f}{\partial W_{\text{output}}} = G_{\text{out}} \dfrac{\partial \hat y}{\partial W_{\text{output}}} \\[2mm]
W_{H_4} &
G_{H_4} = \dfrac{\partial f}{\partial \hat y} \dfrac{\partial \hat y}{\partial H_4} = \dfrac{\partial f}{\partial H_4} &
\dfrac{\partial f}{\partial W_{H_4}} = G_{H_4} \dfrac{\partial H_4}{\partial W_{H_4}} \\[2mm]
W_{H_3} &
G_{H_3} = \dfrac{\partial f}{\partial \hat y} \dfrac{\partial \hat y}{\partial H_4} \dfrac{\partial H_4}{\partial H_3} = \dfrac{\partial L}{\partial H_3} &
\dfrac{\partial f}{\partial W_{H_3}} = G_{H_3} \dfrac{\partial H_3}{\partial W_{H_3}} \\[2mm]
W_{H_2} &
G_{H_2} = \dfrac{\partial f}{\partial \hat y} \dfrac{\partial \hat y}{\partial H_4} \dfrac{\partial H_4}{\partial H_3} \dfrac{\partial H_3}{\partial H_2} = \dfrac{\partial f}{\partial H_2} &
\dfrac{\partial f}{\partial W_{H_2}} = G_{H_2} \dfrac{\partial H_2}{\partial W_{H_2}} \\[2mm]
W_{H_1} &
G_{H_1} = \dfrac{\partial f}{\partial \hat y} \dfrac{\partial \hat y}{\partial H_4} \dfrac{\partial H_4}{\partial H_3} \dfrac{\partial H_3}{\partial H_2} \dfrac{\partial H_2}{\partial H_1} = \dfrac{\partial f}{\partial H_1} &
\dfrac{\partial f}{\partial W_{H_1}} = G_{H_1} \dfrac{\partial H_1}{\partial W_{H_1}}
\end{array}
\]

During local preconditioner iterations, the forward activation from the previous subdomain provides the input needed to evaluate the local derivative $\partial H_i / \partial W_{H_i}$, while the stored backward gradient from the next subdomain (e.g.\ $\partial L/\partial H_{i+1}$) allows us to reconstruct $G_{H_i} = (\partial L/\partial H_{i+1})(\partial H_{i+1}/\partial H_i) = \partial L/\partial H_i$. Thus, each subdomain can form its approximated gradients without re-running the frozen portions of the network.

Because a subdomain cannot evaluate the global network output, unless it contains the output layer, it no longer has a well-defined local objective. We therefore use Objective Function Free Optimization (OFFO) methods\cite{Nocedal-Wright-1999}, i.e.,\ optimizers that require only gradients (SGD, Adam, Adagrad, etc.). As a consequence, we cannot compute $\rho^k$ in Step~10 of Algorithm~\ref{alg:apts}; instead of rejecting steps when the loss increases, we accept them to allow exploration. Only the global TR run adjusts the TR radius at each iteration.

The subdomain optimizer is a constrained version of Adam, which we will refer to as CAdam, whose update never exceeds the subdomain TR radius $\Delta$:
\[
\theta_{t+1} \;=\; \theta_t \;-\; \tilde s_t,
\qquad
\tilde s_t =
\begin{cases}
s_t, & \|s_t\|_{p} \le \Delta,\\[2mm]
\dfrac{\Delta}{\|s_t\|_{2}}\,s_t, & \text{otherwise},
\end{cases}
\]
where $s_t$ is the standard Adam step.  
This imposed constraint and the resulting reduction of the effective step radius guarantee that even over five consecutive iterations, the cumulative update remains within the global trust region.

\section{Numerical Results}
\label{sec:results}
In this section, we present numerical results testing our proposed Inexact and stochastic APTS (IAPTS) method. First, we compare it against Adam optimizer on the MNIST dataset, and then against the SGD optimizer on the CIFAR-10 dataset. All experiments are repeated five times, with different random NN parameter initializations, and are run on the Daint Cray XC50 cluster at the Swiss National Supercomputing Centre (CSCS), using the PyTorch distributed framework with the NVIDIA NCCL backend. Reported plots show the average accuracy and loss over those five runs. The complete source code used to produce the results is available at the project’s GitHub repository: \href{https://github.com/cruzas/DD4ML}{DD4ML}

\subsection{MNIST}\label{sec:MNIST}
For the experiments with the MNIST dataset \cite{MNIST}, we select a batch size of 10\,000 and an NN with four convolutional layers followed by two fully connected layers. The NN is partitioned into $2,4,6$ subdomains, each hosted on a separate GPU.

Adam's learning rate is selected via ten preliminary runs, whereas IAPTS employs standard TR parameters. Specifically, we evaluate each optimizer under the following settings:
\begin{itemize}
  \item \textbf{Adam:} learning rate \(\eta=0.0025\).
  \item \textbf{IAPTS:} to solve the optimization problem on each subdomain, we perform five local CAdam iterations with input learning rate set to $\frac{1}{5}\Delta_G^k$, where $k$ is the IAPTS iteration. IAPTS learning rate is initialized at \(\eta=0.01\) and adaptively bounded between a maximum of 1.0 and a minimum of 0.001 throughout the iterations.
\end{itemize}

\begin{figure}[htbp]
    \centering
    \includegraphics[width=0.6\linewidth]{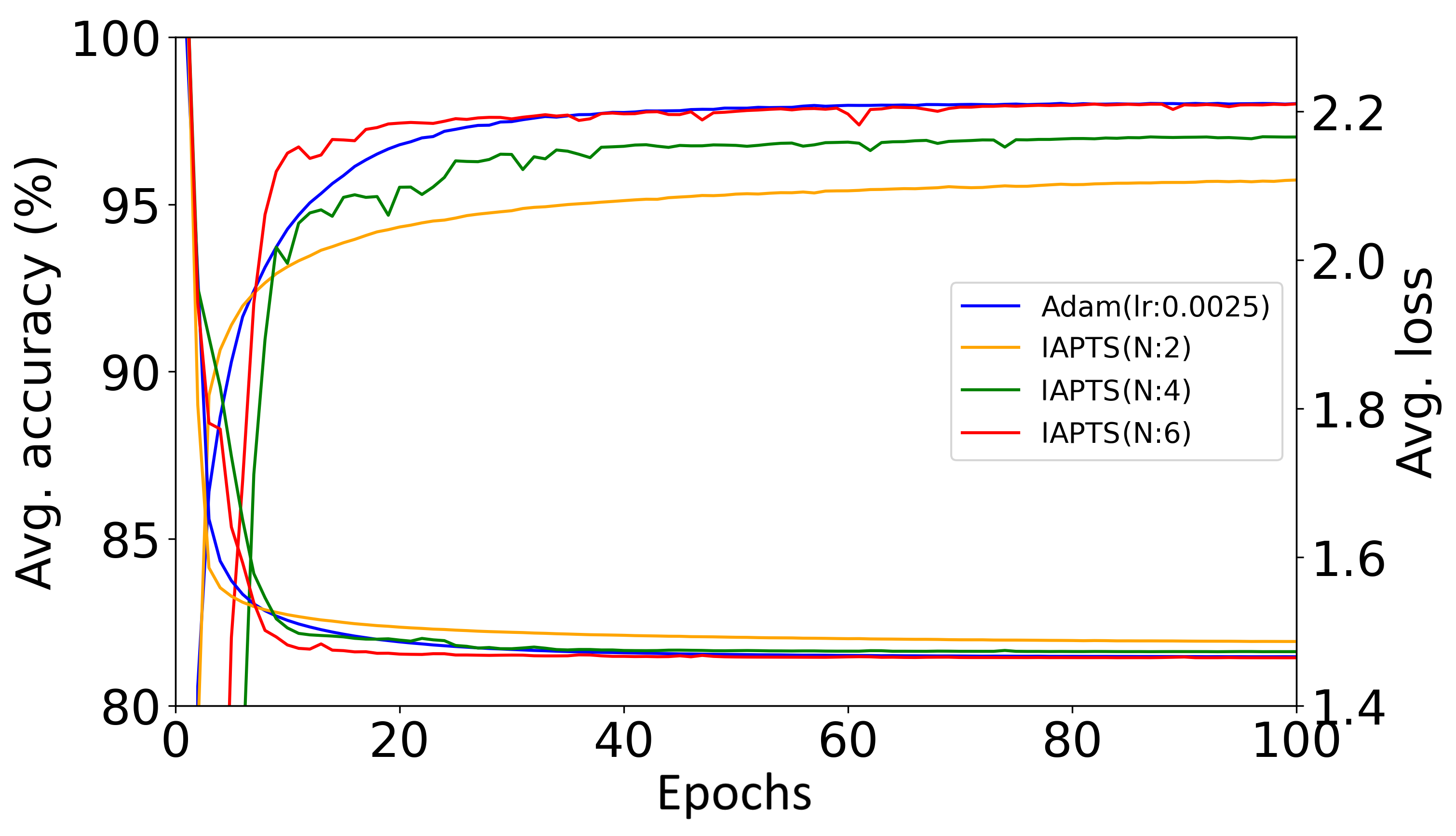}
    \caption{Average training accuracy (left axis) and loss (right axis) over 100 epochs with a batch size of 10\,000 on the MNIST dataset. Solid lines denote the mean across five runs.}
    \label{fig:acc-loss-epochs}
\end{figure}

Figure~\ref{fig:acc-loss-epochs} shows the mean training accuracy (left axis) and loss (right axis) per epoch.
In Figure~\ref{fig:acc-loss-epochs}, IAPTS with 2 subdomains is shown to achieve higher accuracy in the early stages of training, while IAPTS with 4 subdomains surpasses it in the mid and late stages. IAPTS with 6 subdomains attains the highest overall accuracy and remains computationally cheapest per subdomain, suggesting that increasing \(N\) enhances the generalization properties of IAPTS. Similar trends have been observed in related works ~\cite{Trotti-et-al-2025, CruzAlegría-et-al-2025}.

Moreover, IAPTS with 6 subdomains reaches peak accuracy and minimal training loss faster than Adam, although both optimizers converge to comparable final values.

\subsection{CIFAR-10}
Our experiment on the CIFAR‑10 dataset~\cite{cifar10} employs a finely tuned SGD optimizer with learning rate \(\eta = 0.1\) and momentum \(m = 0.9\), a batch size of 200, and the NN is a ResNet with a convolutional layer, followed by 6 residual blocks and finally a fully connected layer. The SGD's hyperparameters were selected based on ten preliminary tuning runs to ensure fast convergence. IAPTS uses the same TR settings as in Section~\ref{sec:MNIST}.
\begin{figure}[htbp]
    \centering
    \includegraphics[width=0.6\linewidth]{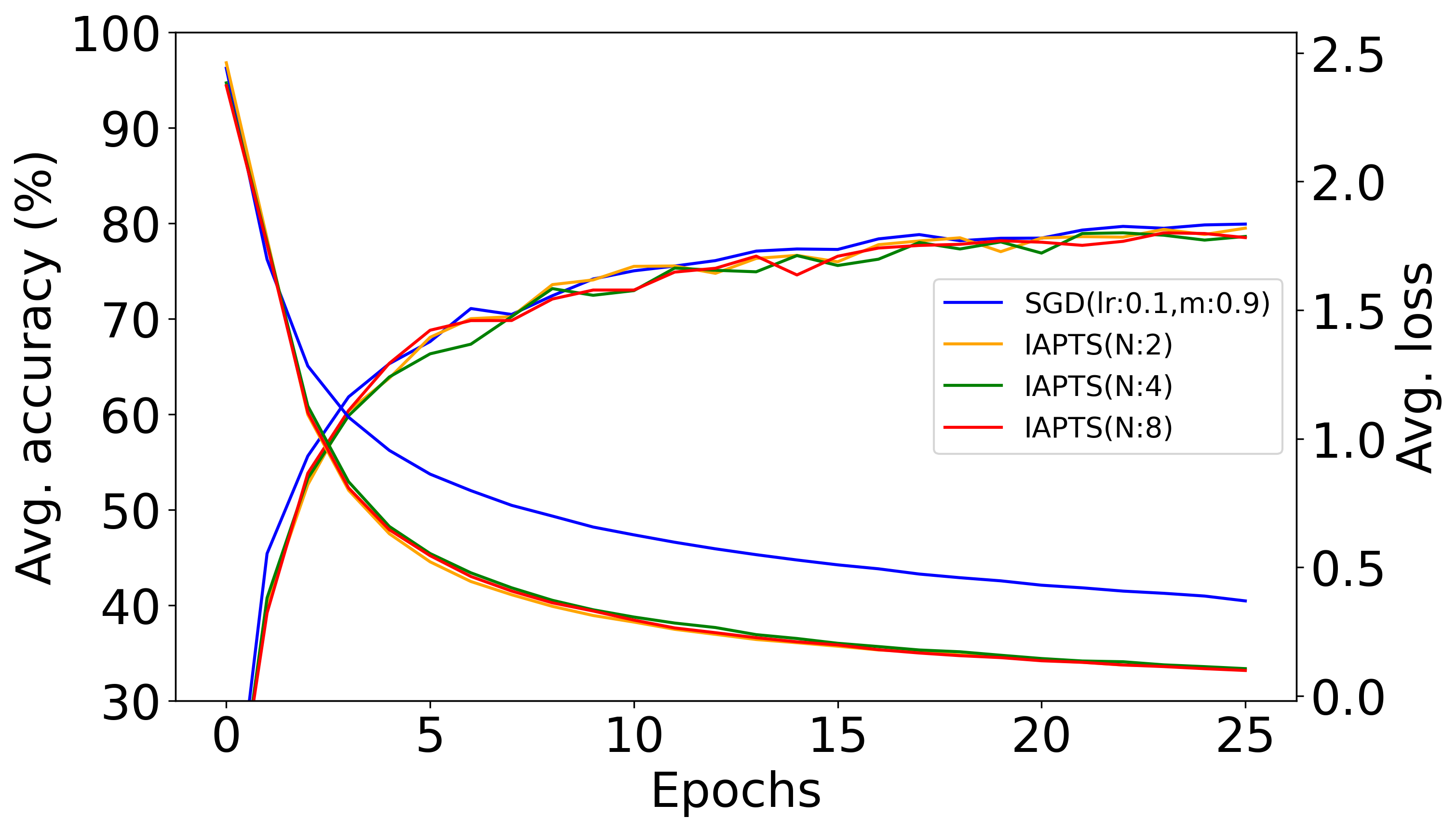}
    \caption{Average training accuracy (left axis) and loss (right axis) over 25 epochs with a batch size of 200 on the CIFAR-10 dataset. Solid lines denote the mean across five runs.}
    \label{fig:acc-loss-epochs2}
\end{figure}
For the CIFAR‑10 dataset, in Figure \ref{fig:acc-loss-epochs2} we note that IAPTS achieves classification accuracies comparable to SGD across all subdomain counts, demonstrating robustness to parallel partitioning. Crucially, from epoch 3 onward, every IAPTS variant maintains a notably lower training loss than SGD, underscoring IAPTS’s accelerated convergence under distributed execution.

It is important to emphasize that IAPTS is not designed to outperform a well-tuned Adam or SGD; rather, it demonstrates that competitive accuracy can be achieved without the need for extensive hyperparameter tuning. Although each IAPTS iteration incurs roughly three times the cost of a single Adam or SGD update (two full forward/backward passes plus multiple inexact backward steps), determining Adam’s and SGD's suitable learning rate required ten preliminary trials. When accounting for the overhead of hyperparameter search, IAPTS remains more efficient overall. Care must still be taken to avoid setting the minimum learning rate too low, as this could slow the convergence of the algorithm.

\section{Conclusion}
\label{sec:conclusion}
Our experiments on MNIST and CIFAR-10 show that IAPTS can reliably attain accuracy levels comparable to those of finely tuned SGD and Adam, without any manual learning rate search. Moreover, increasing the number of subdomains yields a boost in accuracy while decreasing the computational complexity of the algorithm.  Although each IAPTS iteration is more expensive than a single first‑order update, the removal of multiple tuning runs makes IAPTS competitive in overall effort.  In this way, IAPTS delivers a globally convergent\footnote{\label{fn:glob_conv}Global convergence can only be guaranteed in the non-stochastic settings.}, hyperparameter‑robust alternative that scales efficiently in parallel.

\newpage

\bibliographystyle{plain}
\bibliography{pamm-tpl}

\end{document}